\documentclass{aims}
\usepackage{amsmath}
  \usepackage{paralist}
  \usepackage{graphics,dsfont} 
  \usepackage{epsfig} 
 \usepackage[colorlinks=true]{hyperref}
\hypersetup{urlcolor=blue, citecolor=red}

\def\currentvolume{X}
 \def\currentissue{X}
  \def\currentyear{200X}
   \def\currentmonth{XX}
    \def\ppages{X--XX}

\newtheorem{theorem}{Theorem}[section]

\newtheorem*{main}{Main Theorem}
\newtheorem{lemma}[theorem]{Lemma}

\theoremstyle{definition}
\newtheorem{definition}[theorem]{Definition}
\newtheorem{remark}{Remark}

\newcommand{\eps}{\varepsilon}

 \DeclareMathOperator*{\sg} {sg}

\newcommand{\R}{\mathbb{R}}

\newcommand{\seq}[1]{\left\{#1\right\}}
\newcommand{\test}{\varphi}
\newcommand{\Dx}{{\Delta x}}

\newcommand{\norm}[1]{\left\|#1\right\|}
\newcommand{\abs}[1]{\left|#1\right|}
\newcommand{\sgn}{\mathrm{sign}}

\newcommand{\uDx}{u_{\Dx}}

\newcommand{\Pt}{\Pi_T}

\newcommand{\downto}{\downarrow}

\title[Rate of convergence]{An error estimate for viscous approximate solutions to degenerate anisotropic convection-diffusion equations}

\author[ Christian Klingenberg and Ujjwal Koley ]{}


\email{klingenberg@mathematik.uni-wuerzburg.de}
 \email{toujjwal@gmail.com}

\thanks{The second author is supported by Alexander von Humboldt Foundation}

\begin{document}
\maketitle


\centerline{\scshape Christian Klingenberg}
\medskip
{\footnotesize
 \centerline{Institut f\"{u}r Mathematik,}
   \centerline{Julius-Maximilians-Universit\"{a}t W\"{u}rzburg,}
   \centerline{Campus Hubland Nord, Emil-Fischer-Straße 30,}
   \centerline{97074, W\"{u}rzburg, Germany.}
} 

\medskip

\centerline{\scshape Ujjwal Koley}
\medskip
{\footnotesize
 \centerline{Institut f\"{u}r Mathematik,}
   \centerline{Julius-Maximilians-Universit\"{a}t W\"{u}rzburg,}
   \centerline{Campus Hubland Nord, Emil-Fischer-Straße 30,}
   \centerline{97074, W\"{u}rzburg, Germany.}
} 

\bigskip

 \centerline{(Communicated by the associate editor name)}

\begin{abstract}
 We consider a viscous approximation for a nonlinear degenerate convection-diffusion equations in two space dimensions, and prove an $L^1$ error estimate. Precisely, we show that the $L^1_{\mathrm{loc}}$ difference between the approximate solution and the unique entropy solution converges at a rate $\mathcal{O}(\eps^{1/2})$, where $\eps$ is the viscous parameter.
\end{abstract}

\section{Introduction}
\label{sec:intro}
In this paper, we are interested in certain ``viscous'' approximations of entropy solutions of the following Cauchy problem
\begin{equation}
  \label{eq:main}
  \begin{cases}
    u_t + f(u)_y = u_{xx}, &(x,y,t)\in \Pt,\\
    u(x,y,0)=u_0(x,y), & (x,y)\in \R^2,
  \end{cases}
\end{equation}
where $\Pt=\R \times \R \times (0,T)$ with $T>0$ fixed, $u:\Pt\to \R$ is the unknown function and $f:\R \rightarrow \R$ is the convective flux function. The main characteristics of this type of equations is that it has mixed parabolic-hyperbolic type, due to the directional separation of the diffusion and convection effects: while matter is convected along the $y$ axis, it is simultaneously diffused along orthogonal direction. 

The existence of solutions of equation \eqref{eq:main} can be obtained by the classical method of adding a vanishing viscosity, in other terms a diffusion, in the missing direction (along the $y$ axis). Since the most characteristics of Equation \eqref{eq:main} is that it has mixed parabolic-hyperbolic type, or in other words, it has strong degeneracy due to the lack of diffusion in the $x$-direction, it is difficult to establish the uniqueness of solutions of \eqref{eq:main}. In \cite{ez1}, authors have proved the existence as well as uniqueness of such problems. 

In this paper we are interested in certain approximate solutions of \eqref{eq:main} coming from solving the uniformly parabolic problem
\begin{equation}
\begin{aligned}
 u^{\eps}_t + f(u^{\eps})_y =  u^{\eps}_{xx} + \eps  u^{\eps}_{yy}.
\end{aligned}
\label{eq:regular}
\end{equation}
We refer to $u^{\eps}$ as a ``viscous'' approximate solution of \eqref{eq:main}. Since the convergence of $u^{\eps}$ to the unique entropy solution $u$ of \eqref{eq:main} as $\eps \downto 0$ is well known, so our interest here is to give an explicit rate of convergence for $u^{\eps}$ as $\eps \downto 0$, i.e., an $L^1$ error estimate for viscous approximate solutions. There are several ways to prove such an error estimate. One way is to view it as a consequence of a continuous dependence estimate. Combining the ideas of \cite{ez1} with a variant of Kruzkov's ``doubling of variables'' device for \eqref{eq:main}, we prove that $\norm{u^{\eps} -u}_{L^1(\Pt)} = \mathcal{O}(\eps^{1/2})$. Although our proof is of independent interest, it may also shed some light on how to obtain an error estimates for numerical methods.

\section{Preliminaries}
\label{sec:prelim}
Independently of the smoothness of the initial data, due to the lack of diffusion in the $x$-direction, jumps may form in the solution $u$. Therefore we consider solutions in the weak sense, i.e., 
\begin{definition}
  Set $\Pt=\R \times \R \times (0,T)$. 
  A function 
  $$
  u(t,x) \in C\left([0,T];L^1(\R^2)\right) \cap L^\infty(\Pt)
  $$ 
  is a weak solution of the initial value problem 
  \eqref{eq:main} if it satisfies:
  \begin{itemize}
  \item [(a)]For all test functions $\test\in
    \mathcal{D}(\R^2 \times [0,T))$ \label{def:w2}
    \begin{equation}
      \label{eq:weaksol}
      \iiint_{\Pt} \left( u\test_t + f(u) \test_y + u \test_{xx} \right)\,dx \,dy \,dt 
      + \iint_{\R^2} u_0(x,y) \test(x,y,0) \,dx \,dy= 0.
    \end{equation}
   \item [(b)]$u = u(x,y,t)$ is continuous at $t=0$ as a function: $[0,T) \rightarrow L^1(\R^2)$.
  \end{itemize}
\end{definition}

Since weak solutions are not uniquely determined by their initial data, one must impose an additional entropy condition to single out the physically relevant solution. 
\begin{definition}
  A weak solution $u$ of the initial value problem \eqref{eq:main} is 
  called an entropy solution, if the following entropy inequality holds for all test 
  functions $0 \le \test \in \mathcal{D}(\R^2 \times (0,T))$:
\begin{equation}
      \label{eq:entropysol}
 \begin{aligned}
 &\iiint_{\Pt} \abs{u -\psi(x)}\test_t +
   \sgn (u -\psi(x)) (f(u) -f(\psi(x))) \test_y \,dxdydt \\
  &\qquad \ge \iiint_{\Pt} - \abs{u- \psi(x)} \test_{xx}  - \sgn (u- \psi(x)) \psi_{xx} \test \,dxdydt. 
 \end{aligned}     
    \end{equation}
\end{definition}
\begin{remark}
 Note that the above entropy condition is inspired by Kru\v{z}kov's entropy condition for scalar conservation laws. However, due to the presence of diffusion term $u_{xx}$ in \eqref{eq:main}, Kru\v{z}kov's entropy condition has to be modified in this case. This is done by introducing as entropy test functions all functions of the form $\abs{u -\psi(x)}$ and $\psi$ smooth, while Kru\v{z}kov's definition asks for $\psi$ to be a constant.
\end{remark}

Regarding the regularized problem \eqref{eq:regular}, it is well known \cite{ez1} that for all $\eps > 0$ this problem has a unique solution
\begin{align*}
 u^{\eps} \in C([0,T]; L^1(\R^2)) \cap L^{\infty}(\Pt).
\end{align*}
In fact, 
\begin{align*}
 u^{\eps} \in C([0,T]; W^{2,p}(\R^2)) \cap C^{1}([0,T]; L^{p}(\R^2)),
\end{align*}
for all $1 < p < \infty$.
In addition, $u^{\eps}$ satisfies the following properties
\begin{itemize}
 \item  [(a)]$\iint_{\R^2} u^{\eps}(x,y,t) \,dxdy = \iint_{\R^2} u_0(x,y) \,dxdy$,
\item  [(b)]$ \text{TV}\, (u^{\eps}(\cdot,\cdot, t)) \le \text{TV}\, (u_0) $,
\item [(c)]$ \norm{\partial_t u^{\eps}(t)}_1 \le C [\text{TV}\, (u_0) + \text{TV}\, ((u_0)_x)]$.
\end{itemize}

Now we are in a position to state our main result, which is the following
\begin{main}
  Let $u$ be the unique entropy solution to \eqref{eq:main} and $u^{\eps}$
  be as defined by \eqref{eq:regular}. Assume 
  that $u_0 \in BV$, $f$ is Lipschitz continuous. Choose a constant 
  $$
  M>\max_{\abs{u}<\norm{u_0}_{L^\infty(\R^2)}} \abs{f'(u)},
  $$
  and another constant $L>MT$, where $T>0$.
  Then there exists a constant $C$, independent
  of $\eps$, but depending on $f$, $L$, $T$ and $u_0$, such that
  \begin{equation*}
    \int_{\R}\int_{-L+Mt}^{L-Mt}\abs{u(t,x,y) - u^{\eps}(t,x,y)}\,dxdy  \le C \eps^{1/2}
    \quad\text{for $t\le T$.}
  \end{equation*}
\end{main}

\section{Proof of the main theorem}
\label{sec:proof} 

The theorem will be proved by a ``doubling of the variables'' argument, which was introduced by Kruzhkov \cite{kr1,kr2} in the context of hyperbolic conservation laws. 

First, observe that from \eqref{eq:entropysol} we have for any test function $\test$ with compact support in $\R \times \R\times (0,T)$ and $u = u(x,z,s)$,
\begin{equation}
 \begin{aligned}
 &\iiint_{\Pt} \abs{u -\psi(x)}\test_s +
   \sgn (u -\psi(x)) (f(u) -f(\psi(x))) \test_z \,dxdzds \\
  &\qquad \ge \iiint_{\Pt} - \abs{u- \psi(x)} \test_{xx}  - \sgn (u- \psi(x)) \psi_{xx} \test \,dxdzds. 
 \end{aligned}
\label{eq:main_1}
\end{equation}
For the regularized equation \eqref{eq:regular}, we start not with the entropy condition, but in the argument leading up to this condition. To do that, first define the regularized counterpart of the signum function as
$$
\sgn_\eta(\sigma) =
\begin{cases}
  \sgn(\sigma) & \abs{\sigma}>\eta,\\
  \sin\left(\frac{\pi\sigma}{2\eta}\right) &\text{otherwise,}
\end{cases}
$$
where $\eta>0$ and the signum function is defined as
$$
\sgn(\sigma) =
\begin{cases}
  -1 & \sigma<0,\\
  0 & \sigma = 0,\\
  1 & \sigma>0,
\end{cases}
$$

Set
\begin{equation*}
  \psi_\eta(u,\psi)=\int_{\psi}^{u} \sgn_\eta (z-\psi)\,dz.
\end{equation*}
This is a convex entropy for all $\psi$. Set $u^{\eps}=u^{\eps}(x,y,t)$ and
rewrite \eqref{eq:regular} as
\begin{equation*}
  u^{\eps}_t + \left(f(u^{\eps})-f(\psi)\right)_y = u^{\eps}_{xx} + \eps  (u^{\eps} -\psi(x))_{yy},
\end{equation*}
and multiply this with $(\psi_\eta)_u(u,\psi) \, \test$ where $\test$ is a test
function with compact support in $\R \times \R\times (0,T)$.
Observe that the solution $u^{\eps}$ of \eqref{eq:regular} is smooth. Hence after a partial
integration, we arrive at
\begin{align*}
  &\iiint_{\Pt} \psi_\eta(u^{\eps},\psi) \test_t +
   Q_{\eta} (u^{\eps},\psi) \test_y \,dxdydt \\
  &\qquad = \iiint_{\Pt} - \sgn_\eta (u^{\eps}- \psi) u^{\eps}_{xx}\test  - \eps \sgn_\eta (u^{\eps}- \psi) (u^{\eps} -\psi(x))_{yy} \test \,dxdydt,
\end{align*}
where we have used $Q'_{\eta}(u,\psi) = \psi'_\eta(u,\psi) f'(u)$.
Next, taking limit as $\eta \rightarrow 0$, we end up with the parabolic equality
\begin{equation}
\begin{aligned}
  &\iiint_{\Pt} \abs{u^{\eps} -\psi(x)}\test_t +
   \sgn (u^{\eps} -\psi(x)) (f(u^{\eps}) -f(\psi(x))) \test_y \,dxdydt \\
  &\qquad = \iiint_{\Pt} - \sgn (u^{\eps}- \psi(x)) u^{\eps}_{xx} \test  - \eps \, \sgn (u^{\eps}- \psi(x)) (u^{\eps} -\psi(x))_{yy} \test \,dxdydt, \\
& \qquad = \iiint_{\Pt} - \sgn (u^{\eps}- \psi(x)) u^{\eps}_{xx} \test  - \eps \, \abs{u^{\eps} -\psi(x)}_{yy} \test \,dxdydt.
\end{aligned}
\label{eq:main_2}
\end{equation}
At this point we are ready to use ``doubling of the variables'' device.
First, using the entropy inequality \eqref{eq:main_1} for the solution $u=u(x,z,s)$ with $\psi(x)= u^{\eps}(x,y,t)$, we get for $(y,t) \in \R \times (0,T)$
\begin{equation}
 \begin{aligned}
 &\iiint_{\Pt} \abs{u -u^{\eps}}\test_s +
   \sgn (u -u^{\eps}) (f(u) -f(u^{\eps})) \test_z \,dxdzds \\
  &\qquad \ge \iiint_{\Pt} - \abs{u- u^{\eps}} \test_{xx}  - \sgn (u- u^{\eps}) u^{\eps}_{xx} \test \,dxdzds. 
 \end{aligned}
\label{eq:main_3}
\end{equation}
Similarly, from the parabolic equality \eqref{eq:main_2} for the solution $u^{\eps} = u^{\eps}(x,y,t)$ with $\psi(x)= u(x,z,s)$, we get for $(z,s) \in \R \times (0,T)$
\begin{equation}
\begin{aligned}
  &\iiint_{\Pt} \abs{u^{\eps} -u}\test_t +
   \sgn (u^{\eps} -u) (f(u^{\eps}) -f(u)) \test_y \,dxdydt \\
  & \qquad = \iiint_{\Pt} - \sgn (u^{\eps}- u) u^{\eps}_{xx} \test  - \eps \, \abs{u^{\eps} -u}_{yy} \test \,dxdydt.
\end{aligned}
\label{eq:main_4}
\end{equation}
We now integrate \eqref{eq:main_3} over $(y,t) \in \R \times (0,T)$ and \eqref{eq:main_4} over $(z,s) \in \R \times (0,T)$. Addition of those two results yields
\begin{equation}
 \begin{aligned}
 & \int\iiiint_{Q_T} \abs{ u^{\eps} -u } (\test_t +\test_s) +
   \sgn ( u^{\eps} -u ) ( f(u^{\eps}) -f(u)) ( \test_y + \test_z) \,dX \\
  &\qquad \ge \int\iiiint_{Q_T} - \abs{u^{\eps} - u} \test_{xx}  - \eps \, \abs{u^{\eps} -u}_{yy} \test \,dX,\\
& \qquad := \int \iiiint_{Q_T} \mathcal{Q}_1 + \mathcal{Q}_2 \,dX,
 \end{aligned}
\label{eq:main_5}
\end{equation}
where $dX=dx\,dy\,dz\,dt\,ds$ and $Q_T = \R \times \R \times \R \times (0,T) \times (0,T)$.

Following Kruzkov and Kuznetsov \cite{kr1,kr2} we now specify a nonnegative test function $ \test =
\test(x,y,t,z,s)$ defined in $Q_T$. To this end, let $\omega \in
C_{0}^{\infty} (\R)$ be a function satisfying
\begin{equation*}
  \mathrm{supp}(\omega) \subset [-1,1], \qquad
  \omega(\sigma) \ge 0 , \qquad \int_{\R} \omega(\sigma)\, d\sigma = 1,
\end{equation*}
and define $\omega_r(x)=\omega(x/r)/r$. Next, let us choose $\phi \in C_c^{\infty} (\R)$ such that
\begin{equation*}
 \phi =
\begin{cases}
 1, & \abs{x}<1,\\ 0 & \abs{x}\ge 2,
\end{cases}
\end{equation*}
and $ 0  \le \phi \le 1$ when $1 \le \abs{x} \le 2$. Then we define $K_{\beta}(x) = \phi(x/\beta)$. We will let $\beta \rightarrow \infty$ later.  
Furthermore, let $h(z)$ be
defined as
\begin{equation*}
  h(z)=
  \begin{cases}
    0, &z<-1,\\ z+1 & z\in [-1,0],\\ 1 &z>0.
  \end{cases}
\end{equation*}
and set $h_\alpha(z)=h(\alpha z)$. 
Let $\nu<\tau$ be two numbers
in $(0,T)$, for any $\alpha> 0$ define
\begin{equation*}
  \begin{gathered}
    H_{\alpha}(t)= \int_{-\infty}^{t}
    \omega_{\alpha}(\xi)\, d\xi, \\
    \begin{aligned}
      \Psi(y,t) &= \left(H_{\alpha_0} (t - \nu) - H_{\alpha_0} (t -
        \tau)\right) \left(h_\alpha(y-L_l(t)) -
        h_\alpha(y-L_r(t)-\frac{1}{\alpha})\right)\\
      &=:{\chi^{\alpha_0}_{(\nu,\tau)}(t)}\,{\chi^\alpha_{(L_l,L_r)}(y,t)}
    \end{aligned}
  \end{gathered}
\end{equation*}
where the lines $L_{l,r}$ are given by
$$
L_l(t)=-L+Mt,\ L_r(t)=L-Mt
$$
where $M$ and $L$ are positive numbers, $M$ will be specified below.
With $0<r< \min\seq{\nu, T - \tau}$ and
$\alpha_0\in(0,\min\seq{\nu-r,T-\tau-r})$ we set
\begin{equation}
  \label{eq:testfn}
  \test(x,y,t,z,s)= K_{\beta}(x)\, \Psi(y,t)\,\omega_r(y-z)\,\omega_{r_0}(t-s).
\end{equation}
We note that $\test$ has compact support and also that we have,
\begin{align*}
  \test_t + \test_s &= K_{\beta}(x) \, \Psi_t(y,t)\,
  \omega_r(y-z)\,\omega_{r_0}(t-s),\\
  \test_y+\test_z&= K_{\beta}(x)\, \Psi_y(y,t)\, \omega_r(y-z)\,\omega_{r_0}(t-s).
\end{align*}
For the record, we note that
\begin{equation}
  \label{eq:Psider}
  \begin{aligned}
    \Psi_t(y,t)&=-\chi^{\alpha_0}_{(\nu,\tau)}(t) M
    \left(h_\alpha'(y-L_l(t)) + h_\alpha'(y-L_r(t)-\frac{1}{\alpha})\right) \\
    & \qquad + 
    \left(\omega_{\alpha_0}(t-\nu) - \omega_{\alpha_0}(t-\tau)\right)
    \chi^\alpha_{(L_l,L_r)}(y,t),\\
    \Psi_y(y,t)  &= \chi^{\alpha_0}_{(\nu,\tau)}(t) 
    \left(h_\alpha'(y-L_l(t)) - h_\alpha'(y-L_r(t)-\frac{1}{\alpha})\right).
  \end{aligned}
\end{equation}
We shall let all the ``small parameters'' $ \alpha$, $\alpha_0$,  $r$,
$r_0$, $\eps$ and $\Dx$ be sufficiently small and the ``large parameter'' $\beta$ be sufficiently big, but fixed. 

Starting the first term on the left of \eqref{eq:main_5}, we
write
\begin{align*}
  \int_{Q_T} & \abs{u^{\eps}-u}  \left(\test_s  +\test_t\right)\,dX \le
  \underbrace{\int_{\Pt} \abs{u^{\eps}(x,y,t)-u(x,y,t)} K_{\beta}(x)\Psi_t \,dxdydt}_\delta \\
  &\quad + \underbrace{\int_{\Pt}\int_{\R} \abs{u(x,y,t)-u(x,y,s)} K_{\beta}(x) \abs{\Psi_t(x,t)}
  \omega_{r_0}(t-s) \,dxdydsdt}_{\beta} \\ 
  &\quad + \underbrace{\int_{Q_T}  \abs{u(x,y,s)-u(x,z,s)} K_{\beta}(x) \abs{\Psi_t(x,t)}
  \omega_{r_0}(t-s)\,\omega_r(x-y)\, dX}_\gamma . 
\end{align*}
Following \cite{uk1}, it is easy to find that
\begin{equation}
  \beta+\gamma \le C\left(r_0+r\right).\label{eq:betagamma}
\end{equation}
To continue the estimate with the first term on the left of
\eqref{eq:main_5}, we split $\delta$ as follows
\begin{align*}
  \delta &= -\iiint_{\Pt} \chi_{(\nu,\tau)}^{\alpha_0}(t) M
    \left(h'_\alpha(y-L_l(t)) + h'_\alpha(y-L_r(t)-\frac{1}{\alpha})\right) \\
& \qquad \qquad \qquad \qquad K_{\beta}(x) \abs{u^{\eps}(x,y,t)-u(x,y,t)} \, dxdydt \\
 &\quad \quad + \iiint_{\Pt} \chi^\alpha_{(L_l,L_r)}(y,t)
  \abs{u^{\eps}(x,y,t)-u(x,y,t)} \\
& \qquad \qquad \qquad \qquad K_{\beta}(x)
  \left(\omega_{\alpha_0}(t-\nu)-\omega_{\alpha_0}(t-\tau)\right) \,dxdydt, \\
& := \delta_1 + \delta_2.
\end{align*}
The term $\delta_1$ will be balanced against the first order
derivative term on the
left hand side of \eqref{eq:main_5}.  To estimate $\delta_2$
we set $e(x,y,t)=\abs{u^{\eps}(x,y,t)-u(x,y,t)}$ and following \cite{uk1}, we find
\begin{equation}
  \label{eq:delta2est}
  \begin{aligned}
    \delta_2 &\le \iint \chi^\alpha_{(L_l,L_r)}(y,\nu) K_{\beta}(x) \abs{u^{\eps}(x,y,\nu)-u(x,y,\nu)}
    \,dxdy \\ & \qquad - \iint \chi^\alpha_{(L_l,L_r)}(y,\tau) K_{\beta}(x)
    \abs{u^{\eps}(x,y,\tau)-u(x,y,\tau)} \,dxdy +
    C\alpha_0.
  \end{aligned}
\end{equation}
Now we rewrite the ``first derivative term'' on the left hand side of
\eqref{eq:main_5}. Doing this, we get
\begin{align*}
  \int_{Q_T} & K_{\beta}(x) \sgn(u^{\eps} -u)
    \left(f(u^{\eps})-f(u)\right)
    \left(\test_y+\test_z\right) \,dX \\ &= 
    \int_{Q_T} K_{\beta}(x) \sg(x,y,z,t,s)\left(f(u^{\eps}(x,y,t))-f(u(x,y,t))\right) \\
   & \qquad \qquad \qquad \qquad  \Psi_y(y,t)\,\omega_r(y-z)\, \omega_{r_0}(t-s)\,dX\\
    &\quad \quad + 
    \int_{Q_T} K_{\beta}(x) \sg(x,y,z,t,s)\left(f(u(x,y,t))-f(u(x,z,s))\right) \\
 & \qquad \qquad \qquad \qquad \Psi_y(y,t)\,
    \omega_r(y-z)\,\omega_{r_0}(t-s)\,dX \\
    &=:\delta_3 + \delta_4,
\end{align*}
where we have set $\sg(x,y,z,t,s)=\sgn(u^{\eps}(x,y,t)- u(x,z,s))$.
We  proceed as follows
\begin{align*}
  \abs{\delta_4} \le \int_{Q_T} &\abs{f(u(x,y,t))-f(u(x,z,s))}
  \chi^{\alpha_0}_{(\nu,\tau)}(t) K_{\beta}(x) \\
  &\qquad \omega_{r_0}(t-s)\,\omega_r(y-z)
  \left(
    h'_\alpha(y-L_l(t)) + h'_\alpha(y-L_r(t)-\frac{1}{\alpha})\right) \,dX.
\end{align*}
We follow \cite{uk1} to estimate each of these two terms to conclude 
\begin{equation}
  \label{eq:delta4bnd}
  \abs{\delta_4} \le C\left(r_0+r\right).
\end{equation}
Again, choosing $M$ larger than the Lipschitz norm of $f$ implies that 
\begin{equation}
  \label{eq:dl1dl2bnd}
  \delta_1+\delta_3\le 0.
\end{equation}

Collecting all the terms we see that
\begin{equation}
  \label{eq:estimate-middle}
  \begin{aligned}
    \iint K_{\beta}(x) \chi^\alpha_{(L_l,L_r)}&(y,\tau) \abs{u^{\eps}(x,y,\tau)-u(x,y,\tau)}
    \,dxdy \\
    &\le \iint K_{\beta}(x) \chi^\alpha_{(L_l,L_r)}(y,\nu)
    \abs{u^{\eps}(x,y,\nu)-u(x,y,\nu)} \,dxdy\\
    &\qquad + C\left(r_0+ r +
      \alpha_0 + \alpha\right) + \biggl|{\int_{Q_T} \mathcal{Q}_1 + \mathcal{Q}_2 \, dX }\biggr|.
  \end{aligned}
\end{equation}
In order to estimate $\mathcal{Q}_1$ and $\mathcal{Q}_2$, we proceed as follows
\begin{align*}
 \abs{\int_{Q_T}\mathcal{Q}_1} \,dX & = \int_{Q_T} \abs{u^{\eps} - u}\test_{xx} \,dX \\
& = \int_{Q_T} \abs{u^{\eps} - u} K_{\beta}''(x) \Psi(y,t)\,\omega_r(y-z)\,\omega_{r_0}(t-s) \,dX \\
& \le \norm{u^{\eps} +u}_{L^{\infty}} \int_{Q_T} \frac{1}{\beta^2} \phi''(\frac{x}{\beta}) \Psi(y,t)\,\omega_r(y-z)\,\omega_{r_0}(t-s) \,dX \\
& \le \frac{C}{\beta^2} \iiint \phi''(\frac{x}{\beta}) {\chi^{\alpha_0}_{(\nu,\tau)}(t)}\,{\chi^\alpha_{(L_l,L_r)}(y,t)} \,dxdydt \\
& \le \frac{CK}{\beta \alpha},
\end{align*}
where $K := \int_{\R} \abs{\phi''(y)} \,dy$.

Next,
\begin{align*}
& \abs{\int_{Q_T}\mathcal{Q}_2} \,dX  = \eps \int_{Q_T} \abs{u^{\eps} - u}_y \test_{y} \,dX \\
& = \eps \int_{Q_T} \sgn({u^{\eps} - u}) (u^{\eps})_y  K_{\beta}(x)\, \Psi_y(y,t)\,\omega_r(y-z)\,\omega_{r_0}(t-s) \\
& + \eps \int_{Q_T} \sgn({u^{\eps} - u}) (u^{\eps})_y K_{\beta}(x)\, \Psi(y,t)\,\omega_r'(y-z)\,\omega_{r_0}(t-s) \,dX \\
& := \mathcal{Q}_{2,1} + \mathcal{Q}_{2,2}.
\end{align*}
Each of the above terms can be approximated as follows
\begin{align*}
& \abs{\int_{Q_T}\mathcal{Q}_{2,1}} \,dX   = \eps \int_{Q_T} \sgn({u^{\eps} - u}) (u^{\eps})_y K_{\beta}(x)\,\Psi_y(y,t)\,\omega_r(y-z)\,\omega_{r_0}(t-s) \,dX \\
& \le \eps K \int_{\Pt}(u^{\eps})_y (x,y,t) \chi^{\alpha_0}_{(\nu,\tau)}(t) 
    \left(h_\alpha'(y-L_l(t)) - h_\alpha'(y-L_r(t)-\frac{1}{\alpha})\right) \,dxdydt \\
& \le K \eps \alpha \iint \abs{u^{\eps}_x} \,dxdy \\
& \le C K \alpha \eps,
\end{align*}
where $K = \norm{K_{\beta}}_{\infty}$ and $C = \text{TV}\, (u^{\eps})$. Similarly for the other term
\begin{align*}
\abs{\int_{Q_T}\mathcal{Q}_{2,2}} \,dX & = \eps \int_{Q_T} \sgn({u^{\eps} - u}) (u^{\eps})_y  K_{\beta}(x)\, \Psi(y,t)\,\omega_r'(y-z)\,\omega_{r_0}(t-s) \,dX \\
& \qquad \le \frac{K \eps}{r} \int_{\Pt}(u^{\eps})_y (x,y,t) {\chi^{\alpha_0}_{(\nu,\tau)}(t)}\,{\chi^\alpha_{(L_l,L_r)}(y,t)}\,dxdydt \\ 
& \qquad  \le C K \frac{\eps}{r},
\end{align*}
where again $K = \norm{K_{\beta}}_{\infty}$ and $C = \text{TV}\, (u^{\eps})$.

Therefore
\begin{equation}
  \label{eq:qsum}
  \abs{\int_{\Pt^2} \mathcal{Q}_1 + \mathcal{Q}_2 \,dX } \le 
  C \eps + \frac{C}{\beta \alpha } + \frac{C \eps}{r} .
\end{equation}
where $C$ depends on (among other things) $L$ and $T$, but not on
the parameters $\alpha_0$, $\alpha$, $r_0$, $r$, $\beta$ or $\eps$.

Now we have proved the follwoing Lemma:
\begin{lemma}
  \label{lem:estimates} Assume that $u$ and $u^{\eps}$ take values in
  the interval $[-K,K]$ for some positive $K$. Let $M>\max_{v\in
    [-K,K]} \abs{f'(v)}$. Then if $T\ge\tau>\nu\ge 0$ and $L-M \tau >0$, we have
  \begin{equation}
    \label{eq:estimates}
    \begin{aligned}
    \int_{\R}\int_{-L+M\tau}^{L-M\tau} &\abs{u^{\eps}(x,y,\tau)-u(x,y,\tau)}\,dxdy \\&\le 
    \int_{\R^2} \abs{u^{\eps}(x,y,\nu)-u(x,y,\nu)} \,dxdy\\
    &\quad + C\biggl[
    r_0 + r  +\alpha + \eps + \frac{1}{\alpha \beta} + \frac{\eps}{r}\biggr].
    \end{aligned}
  \end{equation}
\end{lemma}
This follows from \eqref{eq:estimate-middle} and \eqref{eq:qsum},
observing that we can send $\alpha_0$ to zero. Now we let $u(x,y,t)$ be the unique entropy solution of \eqref{eq:main}. Also note that since we let $\beta$ tends to $\infty$ so $ \gamma = \frac{1}{\beta}$ is small. We set $\alpha=r=r_0$ and $\eps = \gamma$, and assume
that $\alpha$ is sufficiently small, then
\begin{equation}
  \label{eq:vudx-estimate}
  \int_{\R}\int_{-L+Mt}^{L-Mt} \abs{\uDx(x,t)-v(x,t)}\,dx \,dy\le
  C\left(\alpha+\frac{\eps}{\alpha}\right),
\end{equation}
for some constant $C$ which is independent of the small parameters
$\alpha$, $\eps$. This follows from \eqref{eq:estimates}. Then setting $\eps=\alpha^2$ proves the main theorem.

\section*{Acknowledgments}

The second author acknowledges support from the Alexander von Humboldt Foundation, through a Humboldt Research Fellowship for postdoctoral researchers.

\medskip
Received xxxx 20xx; revised xxxx 20xx.
\medskip

\end{document}